\documentclass{amsart}

\usepackage[normalem]{ulem}

\usepackage[a4paper,top=4cm,bottom=4cm,left=2.3cm,right=2.3cm]{geometry}

\usepackage{lineno}

\usepackage{soul}

\usepackage{amssymb, amsmath, amsthm, amsfonts, euscript, enumerate, color}
\usepackage[hypertexnames=false,
    backref=page,
    pdftex,
    pdfpagemode=UseNone,
    breaklinks=true,
    extension=pdf,
    colorlinks=true,
    linkcolor=blue,
    citecolor=blue,
    urlcolor=blue,
]{hyperref}

\usepackage[dvipsnames]{xcolor}
\usepackage{tikz-cd}

\usepackage{amscd}

\usepackage{comment}

\usepackage{bbm}
\usepackage{color,soul}
\usepackage[shortlabels]{enumitem}
\usepackage[english]{babel}

\usepackage{setspace}

\theoremstyle{plain}
\newtheorem{theorem}{Theorem}[section]
\newtheorem{lemma}[theorem]{Lemma}
\newtheorem{proposition}[theorem]{Proposition}
\newtheorem{corollary}[theorem]{Corollary}

\newtheorem*{theorem*}{Theorem}

\theoremstyle{definition}
\newtheorem{definition}[theorem]{Definition}

\newtheorem{example}[theorem]{Example}
\newtheorem*{example*}{Example}
\newtheorem*{definition*}{Definition}

\theoremstyle{remark}
\newtheorem{remark}[theorem]{Remark}

\newcommand{\R}{\mathbb{R}} % \R   = numeri reali
\newcommand{\C}{\mathbb{C}} % \C   = numeri complessi

\usepackage{graphicx} % Required for inserting images

\title{Positive Hermitian curvature flow on 2-step nilpotent Lie groups}
\author{Ettore Lo Giudice}
\address[Ettore Lo Giudice]{
Dipartimento di Scienze Matematiche, Fisiche e Informatiche
Unit\`a di Matematica e Informatica\\
Universit\`a degli Studi di Parma\\
Parco Area delle Scienze 53/A, 43124\\
Parma, Italy}
\email{ettore.logiudice@unipr.it}

\thanks{The author is partially supported by GNSAGA of INdAM}

\begin{document}

\maketitle

\begin{abstract}
    We study the positive Hermitian curvature flow for left-invariant metrics on $2$-step nilpotent Lie groups with a left-invariant complex structure $J$. We  describe the long-time behavior of the flow under the assumption that $J[\mathfrak{g}, \mathfrak{g}]$ is contained in the center of $\mathfrak{g}$. 
    We show that under our assumption the flow $g_{t}$ exists for all positive $t$ and $(G,(1+t)^{-1}g_{t})$ converges, in the Cheeger-Gromov topology, to a $2$-step nilpotent Lie group with a non flat semi-algebraic soliton. Moreover, we prove that, in our class of Lie groups, there exists at most one semi-algebraic soliton solution, up to homothety. Similar results were proved by M. Pujia and J. Stanfield for nilpotent complex Lie groups \cite{P2021, S2021}. In the last part of the paper we study the Hermitian curvature flow for the same class of Lie groups. 
\end{abstract} 

\section{Introduction}\label{Introduction}
In \cite{ST2011}, J. Streets and G. Tian introduced a new family of geometric flows called {\em Hermitian curvature flows} ({\em HCFs}) which generalize the {\em K\"ahler-Ricci flow} to the Hermitian setting. 

Let $(M,J)$ be a complex manifold, then the evolution of a Hermitian metric $g$ on $M$ under a Hermitian curvature flow is the following
\begin{equation*}
    \partial_t g_{t} = - \big( S(g_{t}) - Q(g_{t}) \big), \quad g_{t}|_{t=0} = g,
\end{equation*}
where $S(g)$ denotes the second Chern-Ricci curvature tensor of $g$ and $Q(g)$ is a $(1,1)$-symmetric tensor which is quadratic in the torsion $T$ of the Chern connection of $g$. We recall that $S(g)$ is the $(1,1)$-symmetric tensor defined by 
\begin{equation*}
    S_{j \overline{k}} = g^{\overline{r}s} \Omega_{s \overline{r} j \overline{k}},
\end{equation*}
where $\Omega$ is the curvature tensor of the Chern connection. 

Since the tensor $Q$ does not affect the parabolicity of the flow, it can be chosen to preserve different geometric properties. Originally, in \cite{ST2011}, the tensor $Q$ was chosen to obtain, in the compact case, a gradient flow stable near K\"ahler-Einstein metrics with non-positive scalar curvature. 

In \cite{ST2010}, the tensor $Q$ is chosen in order to preserve the {\em pluriclosed} condition $\partial \overline{\partial}\omega = 0$, where $\omega$ denotes the fundamental form of the Hermitian metric $g$. In \cite{U2019}, Y. Ustinovskiy considered the {\em positive Hermitian curvature flow} ($HCF_{+}$)
\begin{equation}\label{Positive HCF flow}
    \partial_t g_{t} = - \Theta(g_{t}), \quad g_{t}|_{t=0} = g,
\end{equation}
where  
\begin{equation}\label{Tensor of the flow}
    \Theta(g) \doteq S(g) + \frac{1}{2}Q^{2}(g), \quad Q^{2}_{j \overline{k}} \doteq g^{\overline{p}q}g^{\overline{r}s} T_{sq\overline{k}} T_{\overline{rp}j},
\end{equation}
and $T_{js \overline{p}} \doteq g_{l\overline{p}}T^{l}_{js}$. The coefficients $T^{l}_{js}$ denote the components of the torsion of the Chern connection of $g$. This flow preserves the Griffiths-positivity and the dual Nakano-positivity of the tangent bundle.

One of the main motivation to study the positive Hermitian curvature flow is that it is conformally equivalent to the {\em type IIB flow} introduced by D. H. Phong, S. Picard and Z. Zhang in \cite{PPZ2019}. 
We recall that, if we suppose that $(M,J)$ is equipped also with a holomorphic volume form $\Psi$, the {\em type IIB} flow is the geometric flow of Hermitian metrics  
\begin{equation}\label{flow}
\partial_t(\|\Psi\|_{\omega_{t}}\,\omega_{t}^{n-1})=i\partial \bar\partial \omega_{t}^{n-2}\,, \quad \omega_{t}|_{t=0} = \omega.
\end{equation}
This flow preserves the conformally balanced condition 
$$
d(\|\Psi\|_{\omega_{t}}\omega_{t}^{n-1})=0
$$
and under this assumption its stationary points are Calabi-Yau metrics. 
On complex $3$-folds the flow is a special case of the {\em anomaly flow} \cite{PPZ2018, PPZ2018Strominger}. By \cite{FP2020}, if $\omega_t$ is a conformally balanced solution to the anomaly flow, then 
the Hermitian metric induced by 
$$
\eta_t=\|\Psi\|_{\omega_{t}}\,\omega_{t}
$$ 
solves the $HCF_+$.  

\medskip 
The purpose of the present paper is to study the behavior of the positive Hermitian curvature flow \eqref{Positive HCF flow} on Lie groups equipped with a left-invariant complex structure. In this setting, the flow preserves the property of the initial metric to be left-invariant and reduces to an ODE. Our first main result is the following 

\begin{theorem}\label{First theorem of the introduction}
    Let $(G,J,g_0)$ be a simply-connected, $2$-step nilpotent Lie group equipped with a left-invariant Hermitian structure $(J,g_{0})$. Assume that the Lie algebra  $(\mathfrak{g},\mu)$ of $G$ is such that $J\mu(\mathfrak{g},\mathfrak{g})$ is contained in the center of $\mathfrak{g}$. Then, the $HCF_{+}$ starting from $g_{0}$, has a long-time solution $g_t$ such that $(G,(1+t)^{-1}g_t)$ converges to a non-flat semi-algebraic soliton $(\overline{G}, \overline{h})$ as $t \to \infty$ in the Cheeger-Gromov topology.
\end{theorem}

The theorem is proved by using the {\em bracket flow technique} introduced by J. Lauret in \cite{L2011} to study the Ricci flow on Lie groups.  The technique allows us to regard the flow as a flow in the space of brackets instead of in the space of the inner products. This is a convenient point of view for several reasons (especially for the study of the limits in Cheeger-Gromov topology) and it was adopted in many papers (see e.g. \cite{AL2019, EFV2015, FLS2024, LPV2020, L2015, P2021,S2021} and the references therein). 

We mention that, to prove Theorem \ref{First theorem of the introduction}, we show that $\Theta(g)$ can be viewed as a moment map for the action of $\text{GL}(\xi,J)$, where $\xi$ is the center of the Lie algebra of $(G,J)$, on the vector space of Lie brackets $\mathfrak{N}$ defined in \eqref{vector space of Lie bracket} (see Section \ref{Proof of the first theorem of the introduction}). We recall that, a result similar to Theorem \ref{First theorem of the introduction} was proved in \cite{P2021,S2021} for complex nilpotent Lie groups. 

\medskip

Notice that the technical assumption $J\mu(\mathfrak g,\mathfrak g)$ contained in the center of $\mathfrak{g}$ in Theorem \ref{First theorem of the introduction} is in particular satisfied if $J$ preserves the center of $\mathfrak g$, but it is in fact more general (see e.g. \cite[Example 3.9]{R2009}). Moreover, also the Kodaira-Thurston surface, equipped with the complex structure as in \cite[Example 4.5]{FV2015}, satisfies this assumption. We mention that, in \cite{B2025}, the complex structures that satisfy $J\mu(\mathfrak g,\mathfrak g)$ contained in the center of $\mathfrak{g}$ are called $2$-step nilpotent complex structure, in the sense of \cite{CFGU2000}.

\medskip 
We recall that a sequence of homogeneous manifolds $(M_{k},g_{k})$ converges to a homogeneous manifold $(\overline{M}, \overline{g})$ in the Cheeger-Gromov topology if there exist a sequence $\Omega_{k} \subseteq \overline{M}$ of open neighborhoods of a base point $p \in \overline{M}$ and a sequence of embeddings $\psi_{k}: \Omega_{k} \to M_{k}$ such that $\psi_{k}^{\ast}g_{t_{k}}$ converge to $\overline{g}$ smoothly as $k \to \infty$ and $\Omega_{k}$ eventually contains every compact subset of $\overline{M}$. We also recall that, a Hermitian metric $g$ is called a {\em {\em $HCF_{+}$ soliton}} if it satisfies 
\begin{equation}\label{Soliton for the Hcf in the introduction}
    \Theta(g)= c g + \mathcal{L}_{Z} g,
\end{equation}
where $c \in \R$, $\mathcal{L}$ denotes the Lie derivative, $Z$ is a complete holomorphic vector field and $\Theta$ is defined as in \eqref{Tensor of the flow}. A soliton is called {\em shrinking} if $c>0$, {\em steady} if $c=0$ or {\em expanding} if $c<0$. The $HCF_{+}$ starting from a soliton evolves as 
\begin{equation*}
     g_{t} = k(t) \, \varphi_{t}^{\ast}g,
\end{equation*}
where $k(t) > 0$ and $\varphi_t : G \to G$ are biholomorphisms. Furthermore, if $\varphi_t$ is a family of Lie group automorphisms and $g$ is left-invariant, then we call $g$ a {\em semi-algebraic $HCF_{+}$ soliton}.

 By using \cite{FP2020}, Theorem \ref{First theorem of the introduction} can be applied in order to study the type IIB flow on Lie groups. We have the following 
 \begin{corollary}\label{Second Theorem of the introduction}
    Let $(G,J, \Psi)$ be a simply-connected, $2$-step nilpotent Lie group equipped with a left-invariant complex structure $J$ and a nowhere vanishing, left-invariant, holomorphic $(n,0)$-form $\Psi$. Assume that the Lie algebra $(\mathfrak{g},\mu)$ of $G$ is such that $J\mu(\mathfrak{g},\mathfrak{g})$ is contained in the center of $\mathfrak{g}$, then any left-invariant, balanced solution $\omega_{t}$ to 
    \begin{equation}\label{flow1}
        \partial_t(\|\Psi\|_{\omega_{t}}\,\omega_{t}^{n-1})=i\partial \bar\partial \omega_{t}^{n-2}
    \end{equation}
    is immortal. Moreover $(G,(1+t)^{-1}\omega_{t})$ converges to a non-flat, left-invariant, semi-algebraic soliton $(\overline{G},\overline{\omega})$ in the Cheeger-Gromov topology.
\end{corollary}
In the present paper, we say that $\omega$ is a {\em soliton to the type IIB flow} if the metric induced by $\|\Psi\|_{\omega}\,\omega$ is a $HCF_{+}$ soliton. It is simple to observe that $\omega$ is a soliton to the type IIB equation if and only if 
\begin{equation*}
        \Theta(\|\Psi\|_{\omega} g) = \big(c \|\Psi\|_{\omega} + \mathcal{L}_{Z} (\|\Psi\|_{\omega}) \big) \, g + \|\Psi\|_{\omega} \, \mathcal{L}_{Z} g,
    \end{equation*}
where $c \in \R$, $Z$ is a complete holomorphic vector field and $g$ is the metric induced by $\omega$. We say that a soliton to the type IIB flow is semi-algebraic if the metric induced by $\|\Psi\|_{\omega}\,\omega$ is a semi-algebraic $HCF_{+}$ soliton (see Definition \ref{semi-algebraic soliton to the type IIB flow}).

Note that, by \cite[Theorem 2.7]{BDV2009}, a nilmanifold, i.e., a compact quotient of a connected, simply connected nilpotent Lie group by a lattice, equipped with a left-invariant complex structure, always admits a nowhere vanishing, left-invariant, holomorphic $(n,0)$-form.

\medskip 
The last result of the paper characterizes semi-algebraic $HCF_{+}$ solitons.
\begin{theorem}\label{Third Theorem of the introduction}
    Let $(G,J)$ be a simply-connected, non abelian, $2$-step nilpotent Lie group equipped with a left-invariant complex structure $J$. Assume that the Lie algebra $(\mathfrak{g},\mu)$ of $G$ is such that $J\mu(\mathfrak{g},\mathfrak{g})$ is contained in the center of $\mathfrak{g}$, then, every semi-algebraic $HCF_{+}$ soliton is expanding and unique up to homotheties.
\end{theorem}
Let us mention that by Corollary \ref{semi-algebraic are algebraic} every semi-algebraic soliton to the $HCF_{+}$ is algebraic (see Section \ref{solitons}). Hence, the limit solitons in Theorem \ref{First theorem of the introduction} and Corollary \ref{Second Theorem of the introduction} are in fact algebraic. 

\medskip

This paper is organized as follows. In Section \ref{Preliminaries} we establish the notation used throughout this article and compute the tensor $\Theta$ defined in \eqref{Tensor of the flow} in terms of the structure constants of the Lie algebra. Furthermore, we briefly recall the bracket flow technique. In Section \ref{Proof of the first theorem of the introduction}, we prove Theorem \ref{First theorem of the introduction}. Meanwhile, in Section \ref{Proof of the third Theorem of the introduction} we prove Theorem \ref{Third Theorem of the introduction} and we provide some examples. Finally, in Section \ref{Section 5}, we study the Hermitian curvature flow considered in \cite{ST2011} under the assumption that $J\mu(\mathfrak{g}, \mathfrak{g})$ is contained in the center of $\mathfrak{g}$.

\medskip

{\em \bf{Notation and conventions.}} Throughout this paper, we will adopt the Einstein summation convention for sums over repeated indices unless otherwise stated.

\vskip.3truecm
{\em \underline{Acknowledgments:}} The author would like to sincerely thank Adriano Tomassini for his constant support and encouragement, and Luigi Vezzoni for suggesting the study of this problem during the author’s master’s thesis, for many helpful discussions, and for his continued interest in the paper. Special thanks are also due to Elia Fusi for numerous insightful discussions, for carefully reading the paper, and for providing several comments that improved its clarity and quality. The author is also grateful to James Stanfield for his helpful remarks, which contributed to improving the paper, and to Ramiro Lafuente and Duong Hong Phong for their valuable comments and interest in the paper.

\section{Preliminaries}\label{Preliminaries}
Given a Hermitian manifold $(M,J,g)$, we denote by $T_{\C}M \doteq TM \otimes \C$ its complexified tangent bundle. $T_{\C}M$ has the natural splitting 
 $T_{\C}M = T^{1,0}M \oplus T^{0,1}M$ and the bundle $ \Lambda^{r}_{\C}M$ of complex $r$-forms splits consequently as 
\begin{equation*}
    \Lambda^{r}_{\C}M = \bigoplus_{p+q=r} \Lambda^{p,q}M,
\end{equation*} 
where $\Lambda^{p,q}M \doteq \Lambda^{p}(T^{1,0}M)^{\ast} \otimes \Lambda^{q}(T^{0,1}M)^{\ast}$. The Hermitian metric $g$ specifies a canonical connection, called the {\em Chern connection}, which is the unique affine connection $\nabla$ which preserves both $g$ and $J$ and whose $(1,1)$-part of the torsion tensor $T$ vanishes.
  
\medskip Let us consider a Lie group $G$  with a left-invariant Hermitian structure $(J,g)$. We denote by $(\mathfrak{g},\mu)$ the Lie algebra of $G$. Let $\{Z_{1}, \dots, Z_{n}\}$ be a $g$-unitary, left-invariant $(1,0)$-frame on $G$. The Christoffel symbols of the Chern connection can be easily described in terms of the components of $\mu$. 

From $T(Z_{i},Z_{\overline{j}})=0$ we deduce 
$$   
\nabla_{Z_{i}} Z_{\overline{j}}=\mu(Z_{i}, Z_{\overline{j}})^{1,0}\,,\quad    \nabla_{Z_{\overline{j}}} Z_{i} = \mu(Z_{i}, Z_{\overline{j}})^{0,1}
$$
which imply 
\begin{equation*}
    \Gamma^{\overline{s}}_{i \overline{j}} = \mu^{\overline{s}}_{i \overline{j}}, \quad \Gamma^{s}_{\overline{j} i} = \mu^{s}_{\overline{j} i}.
\end{equation*}
Furthermore, since 
\begin{equation*}
    g(\nabla_{Z_{i}} Z_{j}, Z_{\overline{l}}) = - g(Z_{j}, \nabla_{Z_{i}} Z_{\overline{l}}) = - g \big(Z_{j}, \mu(Z_{i},Z_{\overline{l}}) \big), 
\end{equation*}
we infer 
\begin{equation*}
    \Gamma^{l}_{ij} = -\mu^{\overline{j}}_{i \overline{l}}.
\end{equation*}
Next, we compute the components of the second Chern-Ricci curvature tensor of $g$ in terms of the components of $\mu$. Let $\Omega$ denotes the curvature tensor of the Chern connection $\nabla$. We have that
\begin{equation*}
    \begin{split}
        \Omega_{i \overline{l} j \overline{k}} & = - g(\nabla_{Z_{\overline{l}}} Z_{j}, \nabla_{Z_{i}} Z_{\overline{k}}) + g(\nabla_{Z_{i}} Z_{j},\nabla_{Z_{\overline{l}}} Z_{\overline{k}}) - g(\nabla_{(\mu^{r}_{i \overline{l}} Z_{r} + \mu^{\overline{r}}_{i \overline{l}} Z_{\overline{r}})} Z_{j}, Z_{\overline{k}})  \\
        & = - \mu^{r}_{\overline{l} j} \mu^{\overline{r}}_{i \overline{k}} + \mu^{\overline{j}}_{i \overline{r}} \mu^{k}_{\overline{l} r} + \mu^{r}_{i \overline{l}} \mu^{\overline{j}}_{r \overline{k}} - \mu^{\overline{r}}_{i \overline{l}} \mu^{k}_{\overline{r} j},
    \end{split}
\end{equation*}
and, consequently, 
\begin{equation}\label{Second Chern-Ricci curvature tensor}
    S_{j \overline{k}} = - \mu^{r}_{\overline{s} j} \mu^{\overline{r}}_{s \overline{k}} + \mu^{\overline{j}}_{s \overline{r}} \mu^{k}_{\overline{s} r} + \mu^{r}_{s \overline{s}} \mu^{\overline{j}}_{r \overline{k}} - \mu^{\overline{r}}_{s \overline{s}} \mu^{k}_{\overline{r} j}.
\end{equation}
The $(2,0)$-component of the torsion tensor is given by 
\begin{equation*}
    T(Z_{r}, Z_{s}) = \nabla_{Z_{r}} Z_{s} - \nabla_{Z_{s}} Z_{r} - \mu(Z_{r}, Z_{s}),
\end{equation*}
so, by using the relations between the Christoffel symbols and the structure constants expressed above, we get
\begin{equation*}
    T^{v}_{rs} = - \mu^{\overline{s}}_{r \overline{v}} + \mu^{\overline{r}}_{s \overline{v}} - \mu^{v}_{rs},
\end{equation*}
and by contracting with the metric, we get 
\begin{equation*}
    T_{rs \overline{k}} = - \mu^{\overline{s}}_{r\overline{k}} + \mu^{\overline{r}}_{s \overline{k}} - \mu^{k}_{rs}.
\end{equation*}

\begin{proposition}\label{Proposition 2.6}
    Let $(G,J,g)$ be a $2$-step nilpotent Lie group equipped with a left-invariant Hermitian structure $(J,g)$. Assume that $J\mu(\mathfrak{g},\mathfrak{g})$ is contained in the center of $\mathfrak g$. Then, with respect to a left-invariant, g-unitary $(1,0)$-frame $\{Z_{1}, \dots, Z_{n}\}$ on $G$, we have 
    \begin{equation}\label{HCF positive flow in terms of structure constants}
            \Theta(g)(Z_{j}, Z_{\overline{k}}) = g\Big(\mu(Z_{s},Z_{\overline{r}})^{0,1}, Z_{j}\Big) g\Big(\mu(Z_{\overline{s}},Z_{r})^{1,0}, Z_{\overline{k}}\Big) + \frac{1}{2}g\Big(\mu(Z_{\overline{s}},Z_{\overline{r}} ), Z_{j} \Big)g\Big(\mu(Z_{s},Z_{r}), Z_{\overline{k}} \Big).
    \end{equation}
\end{proposition}
\begin{proof}
    Let $\{Z_{1}, \dots, Z_{n}\}$ be a left-invariant, $g$-unitary $(1,0)$-frame on $G$, from \eqref{Second Chern-Ricci curvature tensor} we have 
    \begin{equation*}
        S_{j \overline{k}} = - \mu^{r}_{\overline{s} j} \mu^{\overline{r}}_{s \overline{k}} + \mu^{\overline{j}}_{s \overline{r}} \mu^{k}_{\overline{s} r} + \mu^{r}_{s \overline{s}} \mu^{\overline{j}}_{r \overline{k}} - \mu^{\overline{r}}_{s \overline{s}} \mu^{k}_{\overline{r} j}.
    \end{equation*}
     Since $G$ is $2$-step nilpotent and $J\mu(\mathfrak{g},\mathfrak{g})$ is contained in the center of $\mathfrak g$, we have that $\mu^{r}_{s \overline{s}} \mu^{\overline{j}}_{r \overline{k}} = 0$. Indeed
     \begin{equation*}
        \mu^{r}_{s \overline{s}} \mu^{\overline{j}}_{r \overline{k}} = g\Big(\mu(\mu^{r}_{s \overline{s}}Z_{r},Z_{\overline{k}}),Z_{j} \Big) =  g\Big(\mu\big(\mu(Z_{s},Z_{\overline{s}})^{1,0}, Z_{\overline{k}} \big),Z_{j} \Big)=0.
    \end{equation*}
    In the same way $\mu^{\overline{r}}_{s \overline{s}} \mu^{k}_{\overline{r} j} = 0$, hence 
    \begin{equation*}
        S_{j \overline{k}} =  - \mu^{r}_{\overline{s} j} \mu^{\overline{r}}_{s \overline{k}} + \mu^{\overline{j}}_{s \overline{r}} \mu^{k}_{\overline{s} r}.
    \end{equation*}
    Furthermore, with respect to $\{Z_{1}, \dots , Z_{n}\}$, we get 
    \begin{equation*}
        \begin{split}
            Q^{2}_{j \overline{k}} & = T_{rs \overline{k}} T_{\overline{rs} j} = (- \mu^{\overline{s}}_{r \overline{k}} + \mu^{\overline{r}}_{s \overline{k}} - \mu^{k}_{rs}) (-\mu^{s}_{\overline{r}j} + \mu^{r}_{\overline{s}j} - \mu^{\overline{j}}_{\overline{rs}}) \\
            & = \mu^{s}_{\overline{r}j} \mu^{\overline{s}}_{r \overline{k}} - \mu^{r}_{\overline{s}j} \mu^{\overline{s}}_{r \overline{k}} + \mu^{\overline{j}}_{\overline{rs}} \mu^{\overline{s}}_{r \overline{k}} - \mu^{s}_{\overline{r}j} \mu^{\overline{r}}_{s \overline{k}} + \mu^{r}_{\overline{s}j} \mu^{\overline{r}}_{s \overline{k}}  - \mu^{\overline{j}}_{\overline{rs}} \mu^{\overline{r}}_{s \overline{k}} + \mu^{k}_{rs} \mu^{s}_{\overline{r}j} - \mu^{k}_{rs} \mu^{r}_{\overline{s}j} + \mu^{k}_{rs} \mu^{\overline{j}}_{\overline{rs}},
        \end{split}
    \end{equation*}
    but, since $G$ is $2$-step nilpotent and $J\mu(\mathfrak{g},\mathfrak{g})$ is contained in the center of $\mathfrak{g}$, we get 
    \begin{equation*}
        Q^{2}_{j \overline{k}} = \mu^{s}_{\overline{r}j} \mu^{\overline{s}}_{r \overline{k}} + \mu^{r}_{\overline{s}j} \mu^{\overline{r}}_{s \overline{k}} + \mu^{k}_{rs} \mu^{\overline{j}}_{\overline{rs}}  = 2 \mu^{s}_{\overline{r}j} \mu^{\overline{s}}_{r \overline{k}} + \mu^{k}_{rs} \mu^{\overline{j}}_{\overline{rs}}.
    \end{equation*}
    Thus,
    \begin{equation*}
        \Theta(g)_{j \overline{k}} = - \mu^{s}_{\overline{r} j} \mu^{\overline{s}}_{r \overline{k}} + \mu^{\overline{j}}_{s \overline{r}} \mu^{k}_{\overline{s} r} + \mu^{s}_{\overline{r}j} \mu^{\overline{s}}_{r \overline{k}} + \frac{1}{2}\mu^{\overline{j}}_{\overline{rs}} \mu^{k}_{rs} = \mu^{\overline{j}}_{s \overline{r}} \mu^{k}_{\overline{s} r} + \frac{1}{2}\mu^{\overline{j}}_{\overline{rs}} \mu^{k}_{rs},
    \end{equation*}
    i.e., 
    \begin{equation*}
            \Theta(g)(Z_{j}, Z_{\overline{k}}) =  \, g\Big(\mu(Z_{s},Z_{\overline{r}})^{0,1}, Z_{j}\Big) g\Big(\mu(Z_{\overline{s}},Z_{r})^{1,0}, Z_{\overline{k}}\Big) + \frac{1}{2}g\Big(\mu(Z_{\overline{s}},Z_{\overline{r}} ), Z_{j} \Big)g\Big(\mu(Z_{s},Z_{r}), Z_{\overline{k}} \Big)\,,
    \end{equation*}
  as required.
\end{proof}
Note that if in the statement of Proposition \ref{Proposition 2.6} we assume that $J$ is abelian, i.e., if $\mathfrak g^{1,0}$ is an abelian Lie algebra (see \cite{ABD2012}), then the condition that $J\mu(\mathfrak{g},\mathfrak{g})$ is contained in the center of $\mathfrak g$ is satisfied since the center of $\mathfrak g$ is $J$-invariant and $\Theta(g)$ reduces to 
\begin{equation}\label{HCF flow for J-abelian complex structure}
        \Theta(g)(Z_{j}, Z_{\overline{k}}) =  g\Big(\mu(Z_{s},Z_{\overline{r}})^{0,1}, Z_{j}\Big) g\Big(\mu(Z_{\overline{s}},Z_{r})^{1,0}, Z_{\overline{k}}\Big).
\end{equation}
\begin{remark}
    We observe that, if we assume that the Lie group $G$ is complex, then the tensor \eqref{HCF positive flow in terms of structure constants} reduces to the one studied in \cite{P2021,S2021}.
\end{remark}

\subsection{The bracket flow technique}\label{The bracket flow technique}
In this section we give a brief outline of the bracket flow approach introduced by J. Lauret in \cite{L2011}
\medskip 

Let $(G,J,g_0)$ be a simply connected Lie group equipped with a left-invariant Hermitian metric $g_{0}$ and a left-invariant complex structure $J$. Let $\mathfrak{g}$ denote the Lie algebra of $G$ and $\mu_0$ the Lie bracket of $\mathfrak{g}$. Since $J$ and $g_{0}$ are left-invariant, their value is determined by $\langle \, , \, \rangle \doteq g_{0}(e)$ and $J|_{\mathfrak{g}}$. 

The Lie bracket $\mu_0$ can be seen as an element of the algebraic variety of Lie bracket 
\begin{equation*}
    \widetilde{\mathfrak{L}} \doteq \{\mu \in \Lambda^2 \mathfrak{g}^{\ast} \, \otimes \,  \mathfrak{g} \, : \enskip \text{$\mu$ satisfies the Jacobi identity and $N_{\mu, J} = 0$} \},
\end{equation*}
where, $N_{\mu,J}$ is the Nijenhuis tensor associated to $\mu$. The space $\widetilde{\mathfrak{L}}$ admits an action of the Lie group $\text{GL}(\mathfrak{g},J) \doteq \{f \in \text{GL}(\mathfrak{g}) \, | \, f \circ J = J \circ f \}$. This action is defined by 
\begin{equation*}
    f \cdot \mu \doteq f \circ \mu(f^{-1} \cdot,f^{-1} \cdot).
\end{equation*}
Let us consider the positive Hermitian curvature flow on $G$ starting from $g_{0}$
\begin{equation}\label{Invariant geometric flow}
    \begin{cases}
        \partial_t g_t = -\Theta(g_t), \\
        g(0) = g_0,
    \end{cases}
\end{equation}
where $\Theta$ is defined as in \eqref{Tensor of the flow}. Since the tensor $\Theta$ is equivariant under biholomorphisms, then, if the flow starts from a left-invariant Hermitian metric $g_0$, it admits a left-invariant solution, i.e., a solution made by left-invariant Hermitian metrics over an interval  $I \subseteq \R$ with $0 \in I$. 

Theorem 1.1 in 	\cite{L2015} guarantees that the left-invariant solution $g_t$ evaluated at the identity can be rewritten in terms of $\langle \, , \, \rangle$ and an element of $\text{GL}(\mathfrak{g},J)$. Specifically, there exists a smooth curve $(f_t)_{t \in I} \in \text{GL}(\mathfrak{g},J)$ such that $f_{0} = \text{Id}_{\mathfrak{g}}$ and
\begin{equation*}
    g_{t}(\cdot, \cdot) = \langle f_t \cdot , f_t \cdot \rangle.
\end{equation*}
Simultaneously,
\begin{equation*}
    \mu_t = f_t \cdot \mu_0 
\end{equation*}
is a solution to the bracket flow 
\begin{equation}\label{Bracket flow}
        \begin{cases}
            \frac{d}{dt}\mu_t = -\pi(\Theta_{\mu_t})\mu_t, \\
            \mu(0) = \mu_0.
        \end{cases}
\end{equation} 
Here, $\pi : \text{End}(\mathfrak{g}) \to \text{End}(\Lambda^2 \mathfrak{g}^{\ast} \otimes \mathfrak{g})$ is defined as  
\begin{equation}\label{Representation pi}
    \pi(E)\mu(\cdot,\cdot) \doteq E\mu(\cdot,\cdot) - \mu(E \cdot, \cdot) - \mu(\cdot,E \cdot) \qquad \forall E \in \text{End}(\mathfrak{g}),
\end{equation}
furthermore, $\Theta_{\mu_t} \in \text{End}(\mathfrak{g})$ is related to the value of $\Theta(g_{t})$ at $e \in G$ by
\begin{equation}\label{Endomorphism Theta that depends on mu}
    \Theta_{\mu_t} \doteq f_t \Theta_{g_t} f_{t}^{-1} , \quad  g_{t}(\Theta_{g_t} \cdot, \cdot) = \Theta(g_t)(\cdot,\cdot).
\end{equation}
In particular, the solutions to \eqref{Invariant geometric flow} and to \eqref{Bracket flow} are defined on the same interval. Moreover, by \cite[Theorem 1.1]{L2015}, the Lie group $(G_{\mu_{t}}, J_{\mu_{t}}, g_{\mu_{t}})$, defined by the data $(\mathfrak{g}, J, \mu_{t}, \langle \cdot, \cdot \rangle)$, is equivariantly, biholomorphically isometric to $(G,J, g_{t})$, where $g_{t}$ is a solution to the $HCF_{+}$ \eqref{Invariant geometric flow} and $\mu_{t}$ is a solution to the bracket flow \eqref{Bracket flow}.

\subsection{Special solutions to the $HCF_{+}$ and to the type IIB flow} \label{solitons}
Important solutions to the $HCF_{+}$ are static metrics and solitons. We recall that a Hermitian metric $g$ on a Lie group $(G,J)$ is called a static metric to the $HCF_{+}$ if 
\begin{equation*}
    \Theta(g)= c g, \quad c \in \R.
\end{equation*}
The corresponding solution to the $HCF_{+}$ evolves only by scaling of the initial metric. 

A left-invariant metric $g$ on a Lie group $(G,J)$ is called a {\em semi-algebraic soliton} to the $HCF_{+}$ if and only if  
\begin{equation}\label{semi-algsoliton}
    \Theta_{g} = c \, \text{Id} + \frac{1}{2} \big(D + D^{t} \big), \quad c \in \R, \; D \in \text{Der}(\mathfrak{g}, \mu), \; [D,J]=0,
\end{equation}
where $(\mathfrak{g},\mu)$ denotes the Lie algebra of $G$ and $D\in \text{Der}(\mathfrak{g},\mu)$ means that $D$ is a derivation of $\mu$.

Moreover, a semi-algebraic soliton is called an {\em algebraic soliton} to the $HCF_{+}$ if 
\begin{equation}\label{algsoliton}
    \Theta_{g} = c \, \text{Id} + D, \quad c \in \R, \; D \in \text{Der}(\mathfrak{g}, \mu), \; [D,J]=0, 
\end{equation}
hence, $D \in \text{Der}(\mathfrak{g}, \mu)$ can be chosen to be self-adjoint.

\medskip
%%%%%%%%%%%%%%%%%%%%%%%%%%%%%% Type IIB flow
By using that the type IIB flow is conformally equivalent to the $HCF_{+}$, we can define solitons to the type IIB flow. Let $(M,J,g,\Psi)$ be a $n$-dimensional complex manifold equipped with a complex structure $J$, a Hermitian metric $g$ whose fundamental form is $\omega$ and a nowhere vanishing holomorphic $(n,0)$-form $\Psi$. Let $\omega_{t}$ be a solution to
\begin{equation}\label{type IIB flow}
\partial_t(\|\Psi\|_{\omega_{t}}\,\omega_{t}^{n-1})=i\partial \bar\partial \omega_{t}^{n-2}, \quad \omega_{0} = \omega
\end{equation}
on $I \subseteq \R$ such that  
\begin{equation*}
    d(\|\Psi\|_{\omega_{t}}\omega_{t}^{n-1})=0, \quad \forall t \in I.
\end{equation*}
Then, from \cite{FP2020}, $\widetilde{g}_{t} \doteq \|\Psi\|_{\omega_{t}} g_{t}$, where $g_{t}$ is the Hermitian metric associated to $\omega_{t}$, evolves according to 
\begin{equation}\label{HCF+ associated to type IIB flow}
    \partial_{t} \widetilde{g}_{t} = - \frac{1}{n-1} \Theta(\widetilde{g}_{t}), 
\end{equation}
where $\Theta(\widetilde{g}_{t})$ is defined as in \eqref{Tensor of the flow}. 

Define $\widetilde{g} \doteq \|\Psi\|_{\omega} \, g$. If $\widetilde{g}$ is a soliton to the $HCF_{+}$, then 
\begin{equation*}
        \Theta(\widetilde{g})= c \widetilde{g} + \mathcal{L}_{Z} \widetilde{g},
\end{equation*}
where $c \in \R$, $\mathcal{L}$ denotes the Lie derivative and $Z$ is a complete holomorphic vector field. Hence,
\begin{equation*}
    \begin{split}
        \Theta(\|\Psi\|_{\omega} g) & = c \|\Psi\|_{\omega} g + \mathcal{L}_{Z} \big(\|\Psi\|_{\omega} g \big)\\ 
        & = c \|\Psi\|_{\omega} g + \mathcal{L}_{Z} \big(\|\Psi\|_{\omega} \big) g + \|\Psi\|_{\omega} \mathcal{L}_{Z} g.
    \end{split}
\end{equation*}
This leads to the following definition.
\begin{definition}\label{soliton to the type IIB flow}
    Let $(M,J,g,\Psi)$ be as above. We say that $g$ is a soliton to the type IIB flow if 
    \begin{equation*}
        \Theta(\|\Psi\|_{\omega} g) = \big(c \|\Psi\|_{\omega} + \mathcal{L}_{Z} (\|\Psi\|_{\omega}) \big) \, g + \|\Psi\|_{\omega} \, \mathcal{L}_{Z} g,
    \end{equation*}
    where $c \in \R$, $\mathcal{L}$ denotes the Lie derivative, $Z$ is a complete holomorphic vector field and $\Theta$ is defined as in \eqref{Tensor of the flow}.
\end{definition}

\begin{definition}\label{semi-algebraic soliton to the type IIB flow}
    Let $(G,J,g,\Psi)$ be a Lie group equipped with a left-invariant Hermitian structure $(J,g)$ and a left-invariant nowhere vanishing holomorphic $(n,0)$-form $\Psi$. We say that $g$ is a semi-algebraic soliton to the type IIB flow if 
    \begin{equation*}
        \Theta(\|\Psi\|_{\omega} g) (\cdot, \cdot) = c \, \|\Psi\|_{\omega} \,g(\cdot, \cdot) + \frac{1}{2} \|\Psi\|_{\omega} \Big( g(D \cdot , \cdot) + g(\cdot, D \cdot) \Big), \quad c \in \R, \; D \in \text{Der}(\mathfrak{g}, \mu), \; [D,J]=0,
    \end{equation*}
    where $(\mathfrak{g},\mu)$ denotes the Lie algebra of $G$ and $D\in \text{Der}(\mathfrak{g},\mu)$.

    If $D$ can be chosen to be self-adjoint, we call $g$ an algebraic soliton to the type IIB flow.
\end{definition}

\section{Proof of Theorem $\ref{First theorem of the introduction}$}\label{Proof of the first theorem of the introduction}

According to the assumptions of Theorem $\ref{First theorem of the introduction}$, let $(G,J,g_{0})$ be a simply connected, $2$-step nilpotent Lie group equipped with a left-invariant Hermitian structure. Let $(\mathfrak{g},\mu_{0})$ denote the Lie algebra of $G$ and suppose that $J\mu_{0}(\mathfrak{g},\mathfrak{g})$ is contained in the center of $\mathfrak{g}$. Let $\{Z_{1}, \dots , Z_{n}\}$ be a $g_{0}$-unitary $(1,0)$-frame on $(\mathfrak{g},\mu_{0})$. From \eqref{HCF positive flow in terms of structure constants} and \eqref{Endomorphism Theta that depends on mu} we get 
\begin{equation*}
     \Theta_{g_{0}} (Z_{j}) = \mu_{s \overline{r}}^{\overline{j}}\mu_{\overline{s} r}^{l}Z_{l} + \frac{1}{2} \mu_{\overline{sr}}^{\overline{j}} \mu_{sr}^{l} Z_{l}, \qquad \Theta_{g_{0}}(Z_{\overline{j}}) = \mu^{j}_{\overline{s} r} \mu_{s \overline{r}}^{\overline{l}}Z_{\overline{l}} + \frac{1}{2} \mu_{sr}^{j} \mu_{\overline{sr}}^{\overline{l}} Z_{\overline{l}}.
\end{equation*}
Let $\xi$ be the center of $(\mathfrak{g},\mu_{0})$. Note that if $X\in \xi^{\perp}\doteq \{Z \in \mathfrak{g} \, : \, g_{0}(Z,Y) = 0, \, \, \forall Y \in \xi\}$, then $\Theta_{g_{0}}(X)=0$. Furthermore, $\Theta_{g_{0}}$ preserves the center $\xi$, i.e., if $X \in \xi$, then $\Theta_{g_{0}}(X) \in \xi$. Thus, with respect to the block representation $\mathfrak{g }= \xi^{\perp} \oplus \xi$, the endomorphism $\Theta_{g_{0}}$ has the following form  
\begin{equation}\label{Theta expressed with respect to the center}
    \Theta_{g_{0}} = \begin{pmatrix}
        0 & 0 \\
        0 & \ast
    \end{pmatrix}.
\end{equation}
It can be easily seen that, since \eqref{Theta expressed with respect to the center} holds, then the $HCF_{+}$ starting from a left-invariant metric on $(G,J)$ preserves the splitting $\mathfrak{g}=\xi \oplus \xi^{\perp}$.

\medskip

Let us consider the following space
\begin{equation}\label{vector space of Lie bracket}
    \mathfrak{N} \doteq \{\mu \in \widetilde{\mathfrak{L}} \, : \, \text{$\mu$ is  2-step nilpotent and $J\mu(\mathfrak{g},\mathfrak{g}) \subseteq \xi$} \},
\end{equation}
and consider the group $\text{GL}(\xi,J) \doteq \{f \in \text{GL}(\xi) \, | \, f \circ J = J \circ f\} \subseteq \text{GL}(\mathfrak{g},J)$ on $\mathfrak{N}$. This group is a subgroup of $\text{GL}(\mathfrak{g},J)$ via the embedding $f \mapsto \begin{pmatrix}
\text{Id} & 0\\
0 & f 
\end{pmatrix}$ and the Lie algebra of $\text{GL}(\xi,J)$ is denoted by $\mathfrak{gl}(\xi,J)$. Moreover, define $\mathfrak{p}(\xi, J) \doteq \mathfrak{gl}(\xi,J) \cap \text{Sym}(\xi, J)$, where $\text{Sym}(\xi, J)$ is the set of endomorphisms that commute with $J$ and that are symmetric with respect to the inner product induced by $g_{0}$ on the Lie algebra.
\begin{remark}\label{Remark3.1}
    We recall that on the Lie algebra $\mathfrak{g}$ we can associate a Hermitian structure that is given by the value at the identity of the Lie group $G$ of $g_{0}$ and $J$. We can also associate a Hermitian product to any tensor product of $\mathfrak{g}$ and every tensor product of $\mathfrak{g}^{\ast}$. 
    
    Moreover, let $g$ be another left-invariant Hermitian metric on $(\mathfrak{g},J)$ such that $g(\cdot,\cdot) = g_{0}(A \cdot, A \cdot)$, where $A \in \text{GL}(\mathfrak{g},J)$ and set, following the approach described in Section \ref{The bracket flow technique}, $\mu \doteq A \cdot \mu_{0}$ and 
    \begin{equation*}
        \Theta_{\mu} = A \Theta_{g}A^{-1}, \quad g(\Theta_{g} \cdot , \cdot ) = \Theta(g)(\cdot , \cdot).
    \end{equation*}
    Then, by extending any operator to the complexified tangent bundle and considering a $g_{0}$-unitary, left-invariant frame $\{Z_{1}, \dots,Z_{n}\}$, we get  
    \begin{equation}\label{Formula di Theta_mu 1}
        g_{0}(\Theta_{\mu} Z_{j}, Z_{\overline{k}}) = g (\Theta_{g}A^{-1} Z_{j},A^{-1} Z_{\overline{k}}) = \Theta(g)(A^{-1}Z_{j},A^{-1}Z_{\overline{k}}) = \Theta(g)(\tilde{Z}_{j},\tilde{Z_{\overline{k}}}),
    \end{equation}
    where $\{\tilde{Z_{1}}=A^{-1}Z_1, \dots, \tilde{Z_{n}}=A^{-1}Z_n\}$ is a unitary frame with respect to $g$. Furthermore, 
    \begin{equation}\label{Formula di Theta_mu 2}
        \begin{split}
            \Theta(g)(\tilde{Z}_{j},\tilde{Z_{\overline{k}}}) & = g\Big(\mu_{0}(\tilde{Z}_{s},\tilde{Z}_{\overline{r}})^{0,1}, \tilde{Z}_{j}\Big) g\Big(\mu_{0}(\tilde{Z}_{\overline{s}},\tilde{Z}_{r})^{1,0}, \tilde{Z}_{\overline{k}}\Big) + \frac{1}{2}g\Big(\mu_{0}(\tilde{Z}_{\overline{s}},\tilde{Z}_{\overline{r}} ), \tilde{Z}_{j} \Big)g\Big(\mu_{0}(\tilde{Z}_{s},\tilde{Z}_{r}), \tilde{Z}_{\overline{k}} \Big) \\
            & = g_{0} \Big(\mu(Z_{s},Z_{\overline{r}})^{0,1}, Z_{j}\Big) g_{0}\Big(\mu(Z_{\overline{s}},Z_{r})^{1,0}, Z_{\overline{k}}\Big) + \frac{1}{2}g_{0}\Big(\mu(Z_{\overline{s}},Z_{\overline{r}} ), Z_{j} \Big) g_{0}\Big(\mu(Z_{s},Z_{r}), Z_{\overline{k}} \Big).
        \end{split}
    \end{equation}
\end{remark}
\begin{lemma}\label{Momentum map Lemma}
    Let 
    \begin{equation}\label{Momentum map}
            M : \enskip \mathfrak{N} \setminus \{0\} \to \mathfrak{p}(\xi,J), \quad
            \mu \mapsto \frac{2}{\|\mu\|^2}\Theta_{\mu},
    \end{equation}
    then $M$ is a moment map for the linear action of $\text{GL}(\xi,J)$ on $\mathfrak{N} \setminus \{0\}$, i.e.,
    \begin{equation*}
        \langle \Theta_{\mu}, E \rangle =  \frac{1}{2}\langle \pi(E)\mu, \mu \rangle, \qquad \forall E \in \mathfrak{p}(\xi,J), \, \, \mu \in \mathfrak{N}\setminus \{0\}.
    \end{equation*}
\end{lemma}
\begin{proof}
    Let us consider $\{Z_1, \dots, Z_n\}$ a unitary $(1,0)$-basis of $\mathfrak{g}$ with respect to $\langle\cdot, \cdot\rangle$. From \eqref{Formula di Theta_mu 1} and \eqref{Formula di Theta_mu 2}, we have that 
    \begin{equation*}
        (\Theta_{\mu})_{j}^{l} = \mu_{s \overline{r}}^{\overline{j}}\mu_{\overline{s} r}^{l}  + \frac{1}{2}\mu_{\overline{sr}}^{\overline{j}} \mu_{sr}^{l}, \quad (\Theta_{\mu})_{\overline{j}}^{\overline{l}} = \mu_{\overline{s} r}^{j}\mu_{s \overline{r}}^{\overline{l}} + \frac{1}{2} \mu_{sr}^{j} \mu_{\overline{sr}}^{\overline{l}},
    \end{equation*}
    where $\mu_{s \overline{r}}^{\overline{j}}\mu_{\overline{s} r}^{l}$, $\mu_{\overline{sr}}^{\overline{j}} \mu_{sr}^{l}$ denote the structure constants of $\mu$ with respect to the $\langle \cdot , \cdot \rangle$- unitary basis $\{Z_{1}, \dots, Z_{n}\}$. Note that, since \eqref{Theta expressed with respect to the center} holds, then $\Theta_{\mu} \in \mathfrak{p}(\xi, J)$.
    Furthermore, since $E \in \mathfrak{p}(\xi,J)$, hence, it preserves the elements of $\mathfrak{g}^{1,0}$ and $\mathfrak{g}^{0,1}$, then
    \begin{equation*}
        \langle E,\Theta_{\mu} \rangle =2 \text{Re} \big(E_{j}^{l} (\Theta_{\mu})_{\overline{j}}^{\overline{l}} \big)  = 2 \text{Re} \{ \mu_{\overline{s} r}^{j}\mu_{s \overline{r}}^{\overline{l}} E^{l}_{j} + \frac{1}{2}\mu_{sr}^{j} \mu_{\overline{sr}}^{\overline{l}} E^{l}_{j} \}.
    \end{equation*}
    By using \eqref{Representation pi} and $E \in \mathfrak{p}(\xi,J)$, we get that 
    \begin{equation}\label{Rapresentetion of the endomorphism that preserve the center}
        \begin{split}
            \langle \pi(E)\mu(\cdot , \cdot ) , \mu(\cdot , \cdot ) \rangle & = \langle E \circ\mu(\cdot , \cdot ) , \mu(\cdot , \cdot ) \rangle - \langle \mu(E \cdot , \cdot ) , \mu(\cdot, \cdot ) \rangle - \langle \mu(\cdot , E \cdot ) , \mu(\cdot , \cdot ) \rangle \\
            & = \langle E \circ\mu(\cdot , \cdot ) , \mu(\cdot , \cdot ) \rangle ,
        \end{split}
    \end{equation}
    where the last equality holds true because $E \in \mathfrak{p}(\xi, J)$.
    
    Suppose that $\eta \in \mathfrak{N}$, then
    \begin{equation*}
            \langle \mu, \eta \rangle = 2 \text{Re}\{\mu^{j}_{sr} \eta^{\overline{j}}_{\overline{sr}}\} + 4 \text{Re} \{\mu^{j}_{s \overline{r}} \eta^{\overline{j}}_{\overline{s}r}\},
    \end{equation*}
    so, from \eqref{Rapresentetion of the endomorphism that preserve the center}, we get  
    \begin{equation*}
            \langle \pi(E) \mu, \mu \rangle =  2\text{Re}\{E_{j}^{l} \mu^{j}_{sr}\mu^{\overline{l}}_{\overline{sr}}\} + 4 \text{Re}\{E_{j}^{l}\mu^{j}_{s \overline{r}}\mu^{\overline{l}}_{\overline{s}r} \}.
    \end{equation*}
    Thus 
    \begin{equation*}
        \langle \pi(E) \mu, \mu \rangle = 2 \langle E,\Theta_{\mu} \rangle\,,
    \end{equation*}
    as required.
\end{proof}

Now we are ready to prove Theorem \ref{First theorem of the introduction}. The proof is analogous to \cite[Theorem 3.3]{AL2019}. For the sake of completeness, we provide the proof in our case.
\begin{proof}[Proof of Theorem $\ref{First theorem of the introduction}$]
    The first part of the proof focuses on the long-time existence of the solution, the second part on the convergence to a non-flat semi-algebraic soliton.
    Since 
    \begin{align*}
        \frac{d}{dt}\|\mu_t \|^2  = 2 \langle\frac{d}{dt}\mu_t,\mu_t\rangle = -2\langle\pi(\Theta_{\mu_t})\mu_t, \mu_t\rangle = - 4 \langle \Theta_{\mu_t} , \Theta_{\mu_t}\rangle = -4 \|\Theta_{\mu_t}\|^2 \leq 0,
    \end{align*}
    then, the bracket flow has a long-time solution according to the standard theory of ordinary differential equations. Given that both the bracket flow's solution and the $HCF_{+}$'s solution are defined over the same interval, we can conclude that the $HCF_{+}$ also admits a long-time solution.

    Let us denote by $\nu_t$ the rescaled solution by a factor of $\|\mu_t\|$, i.e. $\nu_t \doteq \frac{\mu_t}{\|\mu_t\|}$, then $\|\nu_t\|=1$. By \cite[Lemma 2.3]{AL2019} we get that $\nu_t$ is a solution to 
    \begin{equation}\label{normalized bracket flow equation}
        \frac{d}{dt}\nu_t = - \pi(\Theta_{\nu_t} + r_{\nu_t} Id_{\mathfrak{g}})\nu_t,
    \end{equation}
    where $r_{\nu_t} \doteq \langle \pi(\Theta_{\nu})\nu,\nu \rangle = 2 \|\Theta_{\nu}\|^2$. Equation \eqref{normalized bracket flow equation} is called the normalized bracket flow equation. 
    
    By means of \cite[Lemma 7.2]{BL2021}, since \eqref{Momentum map} is a moment map, it turns out that \eqref{normalized bracket flow equation} is the negative gradient flow of the following real-analytic functional
    \begin{equation*}
        \begin{split}
            F : \enskip & \mathfrak{N} \setminus 0 \to \R, \quad \mu \mapsto \frac{\|\Theta_{\mu}\|^2}{\|\mu\|^4} .
        \end{split}
    \end{equation*}
    Since the space of unitary brackets is compact, and $\nu_t$ exists $\forall t \in [0, \infty)$, then there must exist an accumulation point $\overline{\nu}$. Thus, the Theorem of \L{}ojasiewicz (\cite{L1993}) ensures that $\nu_t \to \overline{\nu}$ for $t \to \infty$. Consequently, $\overline{\nu}$ is a stationary point of the flow and $\Theta_{\overline{\nu}}+r_{\overline{\nu}}Id_{\mathfrak{g}}$ is a derivation of $\mathfrak{g}$. Thus, the corresponding metric is a algebraic soliton to the $HCF_{+}$. 
    
    A direct computation yields
    \begin{equation}\label{Trace of the endomorphism}
            \text{Tr}(\Theta_{\mu}) = \frac{1}{2} \neq 0,
    \end{equation}
    hence, the soliton is not flat.
    
    The last part of the proof follows from the fact that the asymptotic behaviour of $\|\mu_t\|$ is $ t^{-\frac{1}{2}}$. Specifically, R. Arroyo and R. Lafuente, in the proof of \cite[Theorem 3.3]{AL2019}, demonstrate that the evolution equation of $\|\mu_t\|^2$ is bounded in both directions by an ODE of the form $\frac{d}{dt}y=\widetilde{c}y^{2}$, with $\widetilde{c}>0$. 
     
    According to \cite[Section 2.1]{L2013}, scaling the metric by a factor of $k>0$ is equivalently reflected in scaling the corresponding bracket by a factor of $k^{-\frac{1}{2}}$. Consequently, the asymptotic behavior of the metrics corresponding to $\frac{\mu_t}{\|\mu_t\|}$, which are solutions to \eqref{normalized bracket flow equation}, is $(1+t)^{-1}g_t$.
     
    Finally, by \cite[Corollary 6.20]{L2012}, we know that the convergence of the bracket flow implies the convergence in Cheeger-Gromov topology for the corresponding family of left-invariant metrics. Thus, $(1+t)^{-1}g_t$ subconverges to a non-flat semi-algebraic soliton.
\end{proof}

\section{Proof of Theorem \ref{Third Theorem of the introduction} and examples}\label{Proof of the third Theorem of the introduction}
This section is devoted to the proof of Theorem \ref{Third Theorem of the introduction} and to the construction of some examples.
\begin{proposition}\label{No static metrics}
    Let $(G,J)$ be a simply-connected, non-abelian, $2$-step nilpotent Lie group equipped with a left-invariant complex structure $J$. Assume that the Lie algebra $(\mathfrak{g},\mu)$ of $G$ is such that $J\mu(\mathfrak{g},\mathfrak{g})$ is contained in the center of $\mathfrak{g}$. Then, there are no left-invariant $HCF_{+}$ static metrics.
\end{proposition}
\begin{proof}
    Let us suppose that $g$ is a left-invariant Hermitian static metric to the $HCF_{+}$. Since $G$ is $2$-step nilpotent, then we can choose a $g$-unitary, $(1,0)$-left-invariant frame $\{Z_{1},\dots,Z_{n}\}$, such that $\langle Z_{1} \rangle \perp \xi \otimes \C$, where $\xi$ denotes the center of $(\mathfrak{g},\mu)$. By Proposition \ref{Proposition 2.6}, we get 
    \begin{equation*}
        \Theta(g)(Z_{1},Z_{\overline{1}}) = 0,
    \end{equation*}
    hence $c=0$. Furthermore, we also obtain that 
    \begin{equation*}
        n c = \text{Tr}_{g} \, \Theta(g) = g^{i \overline{j}} \Theta(g)(Z_{i}, Z_{\overline{j}}) = \mu_{s \overline{r}}^{\overline{j}} \mu^{j}_{\overline{s}r} + \frac{1}{2} \mu^{\overline{j}}_{\overline{sr}} \mu^{j}_{sr} = \frac{1}{4} \|\mu\|^{2}.
    \end{equation*}
    Hence $c>0$, which is absurd.
\end{proof}
The proof of Theorem \ref{Third Theorem of the introduction} is analogous to the proof of \cite[Theorem B]{P2021}. For the sake of completeness, we outline a part of the proof in our case.

\begin{proof}[Proof of Theorem \ref{Third Theorem of the introduction}]
    According to the assumptions of Theorem \ref{Third Theorem of the introduction}, let $(G,J)$ be a simply-connected, non abelian, $2$-step nilpotent Lie group equipped with a left-invariant complex structure $J$. Suppose that $(G,J)$ is equipped with a metric $g$ which is a semi-algebraic soliton to the $HCF_{+}$. Denote by $(\mathfrak{g},\mu)$ the Lie algebra of $G$ and suppose that $J\mu(\mathfrak{g},\mathfrak{g})$ is contained in the center of $\mathfrak{g}$. Moreover, consider the $g$-orthogonal splitting of the Lie algebra $\mathfrak{g} = \xi^{\perp} \oplus \xi$, where $\xi$ denotes the center of $(\mathfrak{g},\mu)$. Then 
    \begin{equation}\label{Algebraic soliton}
        \Theta_{\mu} = c \text{Id}_{\mathfrak{g}} + \frac{1}{2} (D + D^{t}),
    \end{equation}
    where the transpose is calculated with respect to the Hermitian product induced on $\mathfrak{g}$.
    In the first part of the proof we show that $D^{t}$ is a derivation. Indeed, since $D$ is a derivation, we get that, if $x \in \xi$, then $D(x) \in \xi$, indeed
    \begin{equation*}
        D\big(\mu(x,y)\big) = \mu\big(D(x),y\big) + \mu \big(x,D(y)\big), \quad \forall y \in \mathfrak{g}.
    \end{equation*}
    So $\mu\big(x,D(y)\big) = 0$, for all $y \in \mathfrak{g}$, hence
    \begin{equation*}
        D = \begin{pmatrix}
            \ast & 0 \\
            \ast & \ast
        \end{pmatrix}, 
    \end{equation*}
    with respect to the decomposition of the Lie algebra above.
    The morphism $\pi$ defined in \eqref{Representation pi} is a Lie algebra homomorphism such that $\pi(E^{t}) = \pi(E)^{t}$. 
    Thus, from Lemma \ref{Momentum map Lemma}
    \begin{equation*}
        \begin{split}
            2 \text{Tr} \, \Theta_{\mu} [D,D^{t}] & = \langle \pi([D,D^{t}]) \mu, \mu \rangle = \langle [\pi(D),\pi(D^{t})] \mu, \mu \rangle \\
            & = \langle (\pi(D) \pi(D)^{t} - \pi(D)^{t}\pi(D)) \mu, \mu \rangle \\
            & = \|\pi(D)^{t} \mu\|^{2} - \|\pi(D) \mu\|^{2} = \|\pi(D^{t}) \mu\|^{2}.
        \end{split}
    \end{equation*}
    The last equality holds because $D$ is a derivation. 
    Furthermore, from \eqref{Algebraic soliton} 
    \begin{equation*}
        \begin{split}
            \text{Tr} \, \Theta_{\mu} [D,D^{t}] = c \, \text{Tr} \, [D,D^{t}] + \frac{1}{2} \text{Tr} \, D[D,D^{t}] + \frac{1}{2} \text{Tr} \, D[D,D^{t}],
        \end{split}
    \end{equation*}
    thus,
    \begin{equation*}
        \begin{split}
            \text{Tr} \, \Theta_{\mu} [D,D^{t}] = 0,
        \end{split}
    \end{equation*}
    so $D^{t}$ is a derivation. 
    We already know that there cannot be a static invariant metric, which implies that $D + D^{t} \neq 0$. Furthermore, a soliton cannot be shrinking, i.e., $c > 0$, because the solution would develop a finite-time singularity (\cite{L2015}). Hence, $c \leq 0$. 
    
    Furthermore, by means of \eqref{Theta expressed with respect to the center} and \eqref{Algebraic soliton}, we have that
    \begin{equation*}
        D + D^{t} = \begin{pmatrix}
            -c\text{Id}_{\xi^{\perp}} & 0 \\
            0 & (D+D^{t})_{\xi}
        \end{pmatrix}.
    \end{equation*}
    Thus, if we suppose that $c=0$, then $\Theta_{\mu} =\frac{1}{2} (D + D^{t})$, hence 
    \begin{equation*}
        \begin{split}
            \text{Tr} \, (D + D^{t})^{2} & = 2\text{Tr} \, \Theta_{\mu} (D + D^{t}) = 2\langle \pi(D + D^{t})\mu,\mu \rangle \\
            & = 2\langle \pi (D) \mu,\mu \rangle + 2\langle \pi (D^{t}) \mu,\mu \rangle = 0,
        \end{split}
    \end{equation*}
    but $D + D^{t}$ is a symmetric operator such that $\text{Tr} \, (D + D^{t})^{2} = 0 $, hence $D + D^{t} = 0$, which is absurd. Thus, every semi-algebraic soliton is expanding, i.e. $c < 0 $.
    
    The proof of the uniqueness follows the same argument of \cite[Theorem B]{P2021}.    
\end{proof}
From the proof of Theorem \ref{Third Theorem of the introduction}, it follows that
\begin{corollary}\label{semi-algebraic are algebraic}
    Let $g$ be a semi-algebraic soliton to the $HCF_{+}$ on the Lie group $(G,J)$, then $g$ is an algebraic $HCF_{+}$ soliton.
\end{corollary}

In the following example, we study the $HCF_{+}$ starting from a {\em balanced} metric on a $6$-dimensional, $2$-step nilpotent Lie algebra. We recall that a Hermitian metric is called balanced if its fundamental form is coclosed.
\begin{example}\label{Example of J-abelian Hcf flow}
    Let $\mathfrak{g}$ be the $6$-dimensional, $2$-step nilpotent Lie algebra which satisfies the following structure equations 
    \begin{equation*}
        \begin{split}
            & de^{i} = 0 \quad i = 1, \dots, 4, \\
            & de^{5} = e^{13} - e^{24}, \quad de^{6} = e^{14} + e^{23},
        \end{split}
    \end{equation*}
    where $e^{ij} \doteq e^{i} \wedge e^{j}$. 
    
    Let us consider the following complex structure 
    \begin{equation*}
        Je^{1} = -e^{2}, \quad Je^{3} = e^{4}, \quad Je^{5} = e^{6},
    \end{equation*}
    then, $Je_{1} = -e_{2}$, $Je_{2} = e_{1}$, $Je_{3}=e_{4}$, $Je_{4} = -e_{3}$, $Je_{5} = e_{6}$, $Je_{6} = -e_{5}$. Thus, $J$ is an abelian complex structure. 
    
    Let us set 
    \begin{equation*}
        Z_{1} \doteq \frac{1}{\sqrt{2}} (e_{1} - iJe_{1}), \quad Z_{2} \doteq \frac{1}{\sqrt{2}} (e_{3} - i Je_{3}), \quad Z_{3} \doteq \frac{1}{\sqrt{2}} (e_{5} - iJe_{5}),
    \end{equation*}
    then, the only non vanishing Lie bracket are:
    \begin{equation*}
        \mu(Z_{1},Z_{\overline{2}}) = - \sqrt{2} Z_{\overline{3}}, \quad \mu(Z_{\overline{1}},Z_{2}) = -\sqrt{2} Z_{3}.
    \end{equation*}
    Thus, if $\{\zeta^{1},\zeta^{2},\zeta^{3}\}$ is dual to the $(1,0)$-frame $\{Z_{1},Z_{2},Z_{3}\}$, the Lie bracket takes the following expression 
    \begin{equation}\label{bracket for the lie algebra}
        \mu = -\sqrt{2} \, \zeta^{1} \wedge \zeta^{\overline{2}} \otimes Z_{\overline{3}} - \sqrt{2} \, \zeta^{\overline{1}}\wedge \zeta^{2}\otimes Z_{3}.
    \end{equation}
    Let us call a metric $g$ {\em diagonal} if it can be written as 
    \begin{equation}\label{diagonal metric}
       g \doteq a \, \zeta^{1} \odot \zeta^{\overline{1}} + b \, \zeta^{2} \odot \zeta^{\overline{2}} + c \, \zeta^{3} \odot \zeta^{\overline{3}},
    \end{equation}
    where $a,b,c \in \R$ and $a,b,c > 0$. We mention that, it is easy to show that $g$ is a balanced metric. With respect to this metric, the orthonormal frame is 
    \begin{equation*}
        \widetilde{Z}_{1} \doteq \frac{1}{\sqrt{a}}Z_{1}, \quad \widetilde{Z}_{2} \doteq \frac{1}{\sqrt{b}} Z_{2}, \quad \widetilde{Z}_{3} \doteq \frac{1}{\sqrt{c}} Z_{3},
    \end{equation*}
    hence 
    \begin{equation*}
        \mu(\widetilde{Z}_{1},\widetilde{Z}_{\overline{2}}) = - \sqrt{\frac{2c}{ab}} \widetilde{Z}_{\overline{3}}, \quad \mu(\widetilde{Z}_{\overline{1}},\widetilde{Z}_2) = - \sqrt{\frac{2c}{ab}} \widetilde{Z}_{3}.
    \end{equation*}
    Let us study the behavior of the $HCF_{+}$ starting at the diagonal metric
    \begin{equation*}
        g_{0} \doteq \zeta^{1} \odot \zeta^{\overline{1}} + \zeta^{2} \odot \zeta^{\overline{2}} + \zeta^{3} \odot \zeta^{\overline{3}}.
    \end{equation*}
    From \eqref{HCF flow for J-abelian complex structure}, we have that, for an arbitrary diagonal Hermitian metric as in \eqref{diagonal metric} 
    \begin{equation*}
        \Theta(g) = \frac{2c^{2}}{ab} \zeta^{3} \odot \zeta^{\overline{3}}.
    \end{equation*}
    Thus, the flow starting from $g_{0}$ is equivalent to 
    \begin{equation*}
        \begin{cases}
            a^{'} = 0, \\
            b^{'} = 0, \\
            c^{'} = - \frac{2c^{2}}{ab},
        \end{cases}
    \end{equation*}
    with $a(0)= b(0)=c(0)=1$. 
    
    Hence 
    \begin{equation*}
        g_{t} = \zeta^{1} \odot \zeta^{\overline{1}} + \zeta^{2} \odot \zeta^{\overline{2}} + \frac{1}{2t + 1} \zeta^{3} \odot \zeta^{\overline{3}}.
    \end{equation*}
\end{example}
\begin{proposition}
    Let $G$ be the simply connected, 2-step nilpotent Lie group with Lie algebra $\mathfrak{g}$ as in Example \eqref{Example of J-abelian Hcf flow}. Then, every left-invariant Hermitian metric on $G$ is an expanding $HCF_{+}$ algebraic soliton.  
\end{proposition}
\begin{proof}
    Let $g$ be a left-invariant Hermitian metric on $G$. Then, we can always find a left-invariant, $(1,0)$-unitary frame $\{W_{1},W_{2},W_{3}\}$ of $g$ such that 
    \begin{equation*}
        W_1 \in \text{Span} \{Z_{1},Z_{2},Z_{3}\} , \quad  W_2 \in \text{Span} \{Z_{2},Z_{3}\}, \quad  W_3 \in \text{Span} \{Z_{3}\},
    \end{equation*}
    where $\{Z_{1},Z_{2},Z_{3}\}$ is the left-invariant $(1,0)$-frame satisfying \eqref{bracket for the lie algebra}. 
    
    With respect to this new frame, we have that 
    \begin{equation*}
        \begin{split}
            & \mu(W_{i},W_{j}) = \mu(W_{\overline{i}},W_{\overline{j}}) = 0, \quad \forall i,j \in \{1,2,3\}, \\
            & \mu(W_{1}, W_{\overline{3}}) = \mu(W_{2},W_{\overline{3}}) = \mu(W_{3},W_{\overline{3}}) = \mu(W_{2},W_{\overline{2}}) = 0,
        \end{split}
    \end{equation*}
    while 
    \begin{equation*}
        \mu(W_{1},W_{\overline{2}}) = u W_{\overline{3}}, \quad \mu(W_{1}, W_{\overline{1}}) = v W_{3} - \overline{v} W_{\overline{3}},
    \end{equation*}
    with $u,v \in \C$ and $u,v \neq 0$.
    Hence 
    \begin{equation*}
        \mu = u \, \alpha^{1} \wedge \alpha^{\overline{2}} \otimes W_{\overline{3}} + \overline{u} \, \alpha^{\overline{1}} \wedge \alpha^{2} \otimes W_{3} + \alpha^{1} \wedge \alpha^{\overline{1}} \otimes (v W_{3} - \overline{v} W_{\overline{3}}),
    \end{equation*}
    where $\{\alpha^{1}, \alpha^{2}, \alpha^{3}\}$ is dual to the left-invariant, $g$-unitary, $(1,0)$-frame $\{W_{1}, W_{2}, W_{3}\}$. From \eqref{HCF flow for J-abelian complex structure}, we have 
    \begin{equation*}
        \Theta(g) = (\|u\|^{2} + \|v\|^{2}) \alpha^{3} \odot \alpha^{\overline{3}}.
    \end{equation*}
    Hence, the endomorphism $\Theta_{g}$ extended to $\mathfrak{g}^{\C}$ is $\Theta_{g} = \text{Diag}(0,0,\|u\|^{2} + \|v\|^{2}, 0, 0, \|u\|^{2} + \|v\|^{2})$. 
    
    Let $D \doteq \Theta_{g} - c \text{Id}_{\mathfrak{g}^{\C}}$, then $D= \text{Diag}(-c,-c,-c + \|u\|^{2} + \|v\|^{2}, -c, -c, -c + \|u\|^{2} + \|v\|^{2})$. Furthermore, we have that  
    \begin{equation*}
        \begin{split}
            & D\mu(W_{1},W_{\overline{2}}) - \mu(D W_{1},W_{\overline{2}}) - \mu(W_{1},DW_{\overline{2}}) = \\
            & = u D W_{\overline{3}} - u D^{1}_{1} W_{\overline{3}} - u D^{5}_{5} W_{\overline{3}} = \\
            & = u(D^{6}_{6} - D^{1}_{1} - D^{5}_{5})W_{\overline{3}},
        \end{split}
    \end{equation*}
    and 
    \begin{equation*}
        \begin{split}
            & D\mu(W_{1},W_{\overline{1}}) - \mu(D W_{1},W_{\overline{1}}) - \mu(W_{1},DW_{\overline{1}}) = \\
            & = D (v W_{3} - \overline{v} W_{\overline{3}}) - D^{1}_{1} (v W_{3} - \overline{v} W_{\overline{3}}) - D^{4}_{4} (v W_{3} - \overline{v} W_{\overline{3}}) = \\
            & = v (D^{3}_{3} - D^{1}_{1} - D^{4}_{4})W_{3} - \overline{v}(D^{6}_{6} - D^{1}_{1} - D^{4}_{4})W_{\overline{3}}.
        \end{split}
    \end{equation*}
    Hence, $D$ is a derivation if and only if $c = - \|u\|^{2} - \|v\|^{2}$.
\end{proof}

\section{Hermitian curvature flow on $2$-step nilpotent Lie groups}\label{Section 5}
Let $(M,J,g)$ be a Hermitian manifold with Chern connection denoted by $\nabla$. Let $T$ be the torsion tensor of the Chern connection $\nabla$ and $Q^{j}$, for $j=1,\dots, 4$ be the $(1,1)$-symmetric tensors, quadratic in the torsion $T$ of the Chern connection $\nabla$, defined by 
\begin{equation}
    \begin{split}
        & Q^{1}_{j \overline{k}} \doteq g^{\overline{p}q}g^{\overline{r}s} T_{js \overline{p}} T_{\overline{kr}q}, \quad Q^{2}_{j \overline{k}} \doteq g^{\overline{p}q}g^{\overline{r}s} T_{sq\overline{k}} T_{\overline{rp}j},\\
        & Q^{3}_{j \overline{k}} \doteq g^{\overline{p}q}g^{\overline{r}s} T_{js \overline{r}} T_{\overline{kp}q}  , \quad Q^{4}_{j \overline{k}} \doteq \frac{1}{2}g^{\overline{p}q}g^{\overline{r}s} \big( T_{qs\overline{r}} T_{\overline{pk}j} + T_{qj \overline{k}} T_{\overline{pr}s} \big),\\
    \end{split}
\end{equation}
where $T_{js \overline{p}} \doteq g_{l\overline{p}}T^{l}_{js}$ and $T^{l}_{js}$ denote the components of $T$.

The Hermitian curvature flow studied in \cite{ST2011} is defined as 
\begin{equation}\label{HCF flow}
    \partial_{t} g_{t} = - K(g_{t}), \quad g_{t}|_{t=0} = g,
\end{equation}
where 
\begin{equation*}
    K(g) \doteq S(g) - Q(g),
\end{equation*}
where $S(g)$ is the second Chern-Ricci curvature tensor of $g$ and $Q(g)$ is defined as follows
\begin{equation}\label{Tensor Q for original HCF}
    Q(g) \doteq \frac{1}{2}Q^{1}(g) - \frac{1}{4}Q^{2}(g) - \frac{1}{2}Q^{3}(g) + Q^{4}(g).
\end{equation}

\medskip

In \cite{LPV2020}, the authors studied the behavior of \eqref{HCF flow} on complex unimodular Lie group and on unimodular Lie group equipped with an abelian complex structure. Here, we study the behavior of \eqref{HCF flow} on $2$-step nilpotent Lie group such that $J \mu(\mathfrak{g},\mathfrak{g})$ is contained in the center of $\mathfrak{g}$.

\medskip Let us consider a Lie group $G$ equipped with a left-invariant Hermitian structure $(J,g)$. We denote by $(\mathfrak{g},\mu)$ the Lie algebra of $G$. Suppose that the Lie algebra $(\mathfrak{g},\mu)$ of $G$ is such that $J\mu(\mathfrak{g},\mathfrak{g})$ is contained in the center of $\mathfrak{g}$.

\medskip Let $\{Z_{1}, \dots, Z_{n}\}$ be a $g$-unitary, left-invariant frame on $G$. Formulae in Section \ref{Preliminaries} imply that 
\begin{equation*}
    \begin{split}
        & Q^{1}_{j \overline{k}} = \mu^{\overline{r}}_{j \overline{s}} \mu^{r}_{\overline{k}s} + \mu^{\overline{j}}_{r \overline{s}} \mu^{k}_{\overline{r}s} + \mu^{s}_{jr} \mu^{\overline{s}}_{\overline{kr}}, \\
        & Q^{2}_{j \overline{k}} = 2 \mu^{s}_{\overline{r}j} \mu^{\overline{s}}_{r \overline{k}} + \mu^{k}_{rs} \mu^{\overline{j}}_{\overline{rs}}, \\
        & Q^{3}_{j \overline{k}} = \mu^{\overline{j}}_{r \overline{r}} \mu^{k}_{\overline{s} s}, \\
        & 2 Q^{4}_{j \overline{k}} = \mu^{\overline{s}}_{r \overline{r}} \mu^{s}_{\overline{k}j} + \mu^{s}_{\overline{r}r} \mu^{\overline{s}}_{j \overline{k}}.
    \end{split}
\end{equation*}
Hence, 
\begin{equation*}
    \begin{split}
        K_{j \overline{k}} & = - \mu^{r}_{\overline{s} j} \mu^{\overline{r}}_{s \overline{k}} + \mu^{\overline{j}}_{s \overline{r}} \mu^{k}_{\overline{s} r} - \frac{1}{2} \big(\mu^{\overline{r}}_{j \overline{s}} \mu^{r}_{\overline{k}s} + \mu^{\overline{j}}_{r \overline{s}} \mu^{k}_{\overline{r}s} + \mu^{s}_{jr} \mu^{\overline{s}}_{\overline{kr}}\big) \\
        & \; \quad + \frac{1}{4} \big(2 \mu^{s}_{\overline{r}j} \mu^{\overline{s}}_{r \overline{k}} + \mu^{k}_{rs} \mu^{\overline{j}}_{\overline{rs}}\big) + \frac{1}{2} \mu^{\overline{j}}_{r \overline{r}} \mu^{k}_{\overline{s} s} - \frac{1}{2} \big(\mu^{\overline{s}}_{r \overline{r}} \mu^{s}_{\overline{k}j} + \mu^{s}_{\overline{r}r} \mu^{\overline{s}}_{j \overline{k}} \big)
    \end{split} 
\end{equation*}
and a direct computation yields 
\begin{equation}\label{Hermitian curvature tensor}
    \begin{split}
        K_{j \overline{k}} & = \frac{1}{2} \big(\mu^{\overline{j}}_{r \overline{s}} \mu^{k}_{\overline{r}s} - \mu^{\overline{r}}_{j \overline{s}} \mu^{r}_{\overline{k}s} - \mu^{s}_{jr} \mu^{\overline{s}}_{\overline{kr}}\big) + \frac{1}{4} \big( \mu^{k}_{rs} \mu^{\overline{j}}_{\overline{rs}} - 2 \mu^{s}_{\overline{r}j} \mu^{\overline{s}}_{r \overline{k}} \big) + \frac{1}{2} \mu^{\overline{j}}_{r \overline{r}} \mu^{k}_{\overline{s} s} - \frac{1}{2} \big(\mu^{\overline{s}}_{r \overline{r}} \mu^{s}_{\overline{k}j} + \mu^{s}_{\overline{r}r} \mu^{\overline{s}}_{j \overline{k}} \big).
    \end{split} 
\end{equation}

\medskip

\begin{remark}
    We observe that, if we assume that the Lie group $G$ is complex and $2$-step nilpotent, then \eqref{Hermitian curvature tensor} reduces to the one studied in \cite{LPV2020}. Moreover, if the Lie group $G$ is nilpotent and it is equipped with an abelian complex structure, then \eqref{Hermitian curvature tensor} reduces to the one studied in \cite[Section 6]{LPV2020}.
\end{remark}

The following theorem provides a proof of a part of \cite[Theorem 6.2]{LPV2020} for our class of Lie groups. To prove the long-time existence of $\partial_{t} g_{t} = - \text{Ric}^{1,1}(g_{t})$, using our techniques, we need an additional assumption (see Theorem \ref{Long-time existence of Ricci}).

\begin{theorem}
    Let $(G,J)$ be a simply-connected, $2$-step nilpotent Lie group equipped with a left-invariant complex structure $J$. Assume that the Lie algebra  $(\mathfrak{g},\mu)$ of $G$ is such that $J\mu(\mathfrak{g},\mathfrak{g})$ is contained in the center of $\mathfrak{g}$. A left-invariant Hermitian metric $g$ is balanced if and only if the trace of the tensor $K$ coincides with the Riemannian scalar curvature. Moreover, if $g$ is balanced, then the tensor $K$ coincides with the $(1,1)$-component of the Riemannian Ricci tensor.
\end{theorem}
\begin{proof}
    Let us denote by $g$ a left-invariant Hermitian metric on $(G,J)$. Moreover, denote by $\nabla^{LC}$ the Levi-Civita connection of $g$ and by $\Gamma^{j}_{ir}$ the Christoffel symbols of $\nabla^{LC}$. Then, the $(1,1)$-component of the Riemannian Ricci tensor is given by 
    \begin{equation*}
        \begin{split}
            \text{Ric}_{j \overline{k}} & = \Gamma^{s}_{r \overline{r}} \Gamma^{k}_{js} + \Gamma^{\overline{s}}_{r \overline{r}} \Gamma^{k}_{j \overline{s}} - \Gamma^{s}_{j \overline{r}} \Gamma^{k}_{rs} - \Gamma^{\overline{s}}_{j \overline{r}} \Gamma^{k}_{r \overline{s}} - \mu^{r}_{js} \Gamma^{k}_{r \overline{s}} + \Gamma^{s}_{\overline{r}r} \Gamma^{k}_{js} + \Gamma^{\overline{s}}_{\overline{r}r} \Gamma^{k}_{j \overline{s}} - \Gamma^{s}_{jr} \Gamma^{k}_{\overline{r}s} - \mu^{r}_{j \overline{s}} \Gamma^{k}_{rs} - \mu^{\overline{r}}_{j \overline{s}} \Gamma^{k}_{\overline{r}s} \\
            & = \big(\Gamma^{s}_{r \overline{r}} + \Gamma^{s}_{\overline{r}r} \big) \Gamma^{k}_{js} + \big( \Gamma^{\overline{s}}_{r \overline{r}} + \Gamma^{\overline{s}}_{\overline{r}r} \big) \Gamma^{k}_{j \overline{s}} - \big(  \Gamma^{s}_{j \overline{r}} + \mu^{r}_{j \overline{s}} \big) \Gamma^{k}_{rs} - \big(\Gamma^{s}_{jr} + \mu^{\overline{r}}_{j \overline{s}} \big) \Gamma^{k}_{\overline{r}s} - \mu^{r}_{js} \Gamma^{k}_{r \overline{s}} - \Gamma^{\overline{s}}_{j \overline{r}} \Gamma^{k}_{r \overline{s}}.
        \end{split}
    \end{equation*}
    By using the Koszul's formula, we get 
    \begin{equation*}
        \begin{split}
            & \Gamma^{s}_{kr} = \frac{1}{2} \big( \mu^{s}_{kr} - \mu^{\overline{k}}_{r \overline{s}} - \mu^{\overline{r}}_{k\overline{s}} \big), \quad \Gamma^{s}_{\overline{k}r} = \frac{1}{2} \big( \mu^{s}_{\overline{k}r} - \mu^{k}_{r \overline{s}} - \mu^{\overline{r}}_{\overline{kl}} \big), \quad \Gamma^{\overline{s}}_{\overline{k}r} = \frac{1}{2} \big( \mu^{\overline{s}}_{\overline{k}r} - \mu^{\overline{r}}_{\overline{k}s} - \mu^{k}_{rs} \big), \\
            & \Gamma^{s}_{k\overline{r}} = \frac{1}{2} \big( \mu^{s}_{k \overline{r}} - \mu^{r}_{k \overline{s}} - \mu^{\overline{k}}_{\overline{rs}} \big), \quad \Gamma^{\overline{s}}_{k\overline{r}} = \frac{1}{2} \big( \mu^{\overline{s}}_{k \overline{r}} - \mu^{r}_{ks} - \mu^{\overline{k}}_{\overline{r}s} \big).
        \end{split}
    \end{equation*}
    Since $G$ is a $2$-step nilpotent Lie group, it is in particular a unimodular Lie group. We recall that, once a unimodular Lie group is equipped with a left-invariant Hermitian structure, we can read the unimodular condition in terms of a left-invariant unitary frame as 
    \begin{equation*}
        \mu^{r}_{ir} + \mu^{\overline{r}}_{i \overline{r}} = 0, \quad i = 1, \dots, n.
    \end{equation*}
    Hence, we get that $\Gamma^{s}_{r \overline{r}} + \Gamma^{s}_{\overline{r} r} = 0$ and $\Gamma^{\overline{s}}_{r \overline{r}} + \Gamma^{\overline{s}}_{\overline{r}r}  = 0$. Moreover,
    \begin{equation*}
        \Gamma^{s}_{j \overline{r}} + \mu^{r}_{j \overline{s}} = \frac{1}{2} \big( \mu^{s}_{j \overline{r}} + \mu^{r}_{j \overline{s}}  - \mu^{\overline{j}}_{\overline{rs}} \big), \quad \Gamma^{s}_{jr} + \mu^{\overline{r}}_{j \overline{s}} = \frac{1}{2} \big( \mu^{s}_{jr} - \mu^{\overline{j}}_{r \overline{s}}  + \mu^{\overline{r}}_{j \overline{s}} \big). 
    \end{equation*}
    By using that the Lie group is $2$-step nilpotent and $J \mu(\mathfrak{g}, \mathfrak{g})$ is contained in the center of $\mathfrak{g}$, we get 
    \begin{equation*}
        \begin{split}
            \text{Ric}_{j \overline{k}} = & \frac{1}{4} \big( \mu^{s}_{jr} \mu^{\overline{s}}_{\overline{rk}} + \mu^{\overline{j}}_{r\overline{s}} \mu^{k}_{\overline{r}s} - \mu^{\overline{r}}_{\overline{s}j} \mu^{r}_{s \overline{k}} - \mu^{s}_{j \overline{r}} \mu^{\overline{s}}_{\overline{k}r} + \mu^{r}_{j \overline{s}} \mu^{\overline{r}}_{s \overline{k}} + \mu^{\overline{j}}_{\overline{rs}} \mu^{k}_{rs} + \mu^{\overline{s}}_{j \overline{r}} \mu^{s}_{r \overline{k}} - \mu^{r}_{js} \mu^{\overline{r}}_{\overline{sk}} + \mu^{\overline{j}}_{\overline{r}s} \mu^{k}_{r \overline{s}} \big) + \frac{1}{2} \mu^{r}_{js} \mu^{\overline{r}}_{\overline{sk}},
        \end{split}
    \end{equation*}
    and a direct computation yields 
    \begin{equation*}
        \text{Ric}_{j \overline{k}} = \frac{1}{2} \big( \mu^{r}_{js} \mu^{\overline{r}}_{\overline{sk}} + \mu^{\overline{j}}_{r\overline{s}} \mu^{k}_{\overline{r}s} - \mu^{\overline{r}}_{\overline{s}j} \mu^{r}_{s \overline{k}} - \mu^{s}_{j \overline{r}} \mu^{\overline{s}}_{\overline{k}r} \big) + \frac{1}{4} \mu^{\overline{j}}_{\overline{rs}} \mu^{k}_{rs}.
    \end{equation*}
    Therefore, 
    \begin{equation*}
        K_{j \overline{k}} - \text{Ric}_{j \overline{k}} = \frac{1}{2} \big(\mu^{\overline{j}}_{r \overline{r}} \mu^{k}_{\overline{s}s} - \mu^{\overline{s}}_{r \overline{r}} \mu^{s}_{\overline{k}j} - \mu^{s}_{\overline{r}r} \mu^{\overline{s}}_{j \overline{k}} \big),
    \end{equation*}
    and 
    \begin{equation*}
        k - \text{Tr}_{g} \text{Ric} = - \frac{1}{2} \mu^{s}_{\overline{r}r} \mu^{\overline{s}}_{l \overline{l}}.
    \end{equation*}
    Since $G$ is unimodular then the metric $g$ is balanced if and only if $\sum_{l} \mu(Z_{l}, Z_{\overline{l}}) = 0$. Therefore, the metric $g$ is balanced if and only if $k = \text{Tr}_{g} \text{Ric}$. 

    Moreover, if $g$ is balanced, then the tensor $K$ coincides with the $(1,1)$-component of the Riemannian Ricci tensor, as required.
\end{proof}

\begin{remark}
    E. Fusi pointed out to me that, by using \cite[Equation 2.8]{AI2001}, one can proved that, in general, for nilpotent Lie groups, a left-invariant Hermitian metric $g$ is balanced if and only if the trace of the tensor $K$ coincides with the Riemannian scalar curvature. 
\end{remark}

Our last result concerns the long-time existence of the solution of $\partial_{t} g_{t} = - \text{Ric}^{1,1}(g_{t})$ for simply-connected, $2$-step nilpotent Lie group $G$ equipped with a left-invariant complex structure $J$ that preserves the commutators of the Lie algebra $(\mathfrak{g},\mu_{0})$ of $G$. We remark that, since the Lie group is $2$-step nilpotent, then $J \mu_{0}(\mathfrak{g}, \mathfrak{g})$ is contained in the center of $\mathfrak{g}$.

\medskip

Note that, if $X \in \mu_{0}(\mathfrak{g}, \mathfrak{g})$, then $P_{g}(X) \in \mu_{0}(\mathfrak{g},\mathfrak{g})$, where $P_{g}$ is the endomorphism associated to $\text{Ric}^{1,1}(g)$ as in \eqref{Endomorphism Theta that depends on mu}. Thus, with respect to the block representation $\mathfrak{g }= \mu_{0}(\mathfrak{g}, \mathfrak{g})^{\perp} \oplus \mu_{0}(\mathfrak{g}, \mathfrak{g})$, the endomorphism $P_{g}$ has the following form  
\begin{equation}\label{P expressed with respect to the center}
    P_{g} = \begin{pmatrix}
        \ast & 0 \\
        \ast & \ast
    \end{pmatrix}.
\end{equation}

Let us consider the following space 
\begin{equation*}
    \mathfrak{V} \doteq \{\mu \in \widetilde{\mathfrak{L}} \, : \, \text{$\mu$ is  2-step nilpotent and $J \mu(\mathfrak{g},\mathfrak{g}) \subseteq \mu_{0}(\mathfrak{g},\mathfrak{g})$ } \},
\end{equation*}
and consider the action of $N \doteq \{ f \in \text{GL} (\mathfrak{g},J) \, | \, f \, \mu_{0}(\mathfrak{g}, \mathfrak{g}) \subseteq \mu_{0}(\mathfrak{g}, \mathfrak{g})\}$ on $\mathfrak{V}$. 

\begin{theorem}\label{Long-time existence of Ricci}
    Let $(G,J)$ be a simply-connected, $2$-step nilpotent Lie group equipped with a left-invariant complex structure $J$. Assume that the Lie algebra  $(\mathfrak{g},\mu)$ of $G$ is such that $J$ preserves the commutators. The parabolic flow $\partial_{t} g_{t} = - \text{Ric}^{1,1}(g_{t})$ has always a long-time solution for every left-invariant initial Hermitian metric.
\end{theorem}
\begin{proof}
    To prove the theorem we use the bracket flow technique, as was done in the proof of Theorem \ref{First theorem of the introduction}. The bracket flow equation associated to this flow is the following 
    \begin{equation*}
        \frac{d}{dt} \mu_{t} = - \pi(P_{\mu_{t}}) \mu_{t}, \quad \mu_{|t=0} = \mu_{0},
    \end{equation*}
    where, for $\mu \in \mathfrak{V}$, we have 
    \begin{equation*}
        (P_{\mu})_{j}^{l} = \frac{1}{2} \big( \mu^{r}_{js} \mu^{\overline{r}}_{\overline{sl}} + \mu^{\overline{j}}_{r\overline{s}} \mu^{l}_{\overline{r}s} - \mu^{\overline{r}}_{\overline{s}j} \mu^{r}_{s \overline{l}} - \mu^{s}_{j \overline{r}} \mu^{\overline{s}}_{\overline{l}r} \big) + \frac{1}{4} \mu^{\overline{j}}_{\overline{rs}} \mu^{l}_{rs}.
    \end{equation*}
    Let us consider a real endomorphism $E$ that commutes with $J$ and such that $E \, \mu_{0}(\mathfrak{g}, \mathfrak{g}) \subseteq \mu_{0}(\mathfrak{g}, \mathfrak{g})$, then 
    \begin{equation*}
        \langle E, P_{\mu} \rangle = 2 \text{Re} \big(E_{j}^{l} (P_{\mu})_{\overline{j}}^{\overline{l}} \big) =  \text{Re} \{E_{j}^{l} \big( \mu^{j}_{\overline{r}s} \mu^{\overline{l}}_{r \overline{s}} - \mu^{r}_{s \overline{j}} \mu^{\overline{r}}_{\overline{s}l} - \mu^{\overline{s}}_{\overline{j}r} \mu^{s}_{l \overline{r}} + \mu^{\overline{r}}_{\overline{js}} \mu^{r}_{sl} \big) \} + \frac{1}{2} \text{Re} \{ E_{j}^{l} \mu^{j}_{rs} \mu^{\overline{l}}_{\overline{rs}} \}.
    \end{equation*}
    Moreover, if $\mu \in \mathfrak{V}$, then
    \begin{equation*}
        \begin{split}
            \langle \pi(E)\mu(\cdot , \cdot) , \mu(\cdot, \cdot) \rangle & = \langle E \circ\mu(\cdot , \cdot ) , \mu(\cdot , \cdot ) \rangle - \langle \mu(E \cdot , \cdot ) , \mu(\cdot, \cdot ) \rangle - \langle \mu(\cdot , E \cdot ) , \mu(\cdot , \cdot ) \rangle \\
            & = 2 \text{Re} \{E_{j}^{l} \mu^{j}_{sr} \mu^{\overline{l}}_{\overline{sr}} \} + 4 \text{Re} \{E_{j}^{l} \big( \mu^{j}_{\overline{r}s} \mu^{\overline{l}}_{r \overline{s}} - \mu^{r}_{sl} \mu^{\overline{r}}_{\overline{sj}} - \mu^{s}_{l \overline{r}} \mu^{\overline{s}}_{\overline{j}r} - \mu^{\overline{r}}_{\overline{s}l} \mu^{r}_{s \overline{j}} \big)\},
        \end{split}
    \end{equation*}
    hence 
    \begin{equation*}
        \langle E, P_{\mu} \rangle = \frac{1}{4} \langle \pi(E)\mu(\cdot , \cdot) , \mu(\cdot, \cdot) \rangle.
    \end{equation*}
    Thus,
    \begin{equation*}
        \frac{d}{dt} \| \mu_{t} \|^{2} = 2 \langle \frac{d}{dt} \mu_{t}, \mu_{t} \rangle = - 2 \langle \pi(P_{\mu_{t}}) \mu_{t} , \mu_{t} \rangle = -8 \langle P_{\mu_{t}} , P_{\mu_{t}} \rangle = - 8 \|P_{\mu_{t}}\|^{2} \leq 0,
    \end{equation*}
    and the thesis follows. 
\end{proof}

\begin{remark}
    We recall that the pluriclosed flow evolves the fundamental form $\omega$ of a Strong K\"ahler with torsion metric $g$ in the direction of the $(1,1)$-component of the Ricci form associated with the Bismut connection of $g$. Meanwhile, in the setting considered in this section, the flow evolves a Hermitian metric in the direction of the $(1,1)$-component of the Ricci tensor of the Levi-Civita connection. We would like to point out that, in \cite[Theorem 1.1]{EFV2015} and \cite[Theorem A]{AL2019}, the authors establish long-time existence and describe the asymptotic behavior of the solution to the pluriclosed flow on $2$-step nilpotent Lie groups endowed with a left-invariant complex structure, without requiring any additional assumptions on the complex structure. Whereas in our case, using our techniques, we need extra assumptions on the complex structure.
\end{remark}


\nocite{*}
\begin{thebibliography}{99}
    \bibitem{ABD2012} \textsc{\small A. Andrada, M. L. Barberis, I. G. Dotti}, \textup{\footnotesize Abelian Hermitian geometry}, \textit{Differential Geom. Appl. }, \textbf{30} (2012), pp. 509--519.

    \bibitem{AI2001} \textsc{\small B. Alexandrov, S. Ivanov}, \textup{\footnotesize Vanishing theorems on Hermitian manifolds}, \textit{Differential Geom. Appl.}, \textbf{14} (2001), pp. 251--265.
    
    \bibitem{AL2019} \textsc{\small R. Arroyo, R. A. Lafuente}, \textup{\footnotesize The long-time behavior of the homogeneous pluriclosed flow}, \textit{Proc. Lond. Math. Soc. (3)}, \textbf{119} (2019), pp. 266--289.

    \bibitem{B2025} \textsc{\small M. L. Barberis}, \textup{\footnotesize Complex structures on two-step nilpotent Lie groups}, \textit{preprint arXiv:2508.05885v1}.
    
    \bibitem{BDV2009} \textsc{\small M. L. Barberis, I. Dotti, M. Verbitsky}, \textup{\footnotesize Canonical bundles of complex nilmanifolds, with applications to hypercomplex geometry}, \textit{Math. Res. Lett.}, \textbf{16} (2009), pp. 331--347.

    \bibitem{BL2021} \textsc{\small C. B\"{o}hm, R. A. Lafuente}, \textup{\footnotesize Real geometric invariant theory}, \textit{London Math. Soc. Lecture Note Ser.}, \textbf{463} (2021), pp. 11--49.

    \bibitem{CFGU2000} \textsc{\small L. A. Cordero, M. Fern\'{a}ndez, A. Gray, L. Ugarte}, \textup{\footnotesize Compact nilmanifolds with nilpotent complex structures: Dolbeault cohomology}, \textit{Trans. Amer. Math. Soc.}, \textbf{352} (2000), pp. 5405--5433.

    \bibitem{EFV2015} \textsc{\small N. Enrietti, A. Fino, L. Vezzoni}, \textup{\footnotesize The pluriclosed flow on nilmanifolds and tamed symplectic forms}, \textit{J. Geom. Anal.}, \textbf{25} (2015), pp. 883--909.

    \bibitem{FP2020} \textsc{\small T. Fei, D. H. Phong}, \textup{\footnotesize Unification of the {K}ahler-{R}icci and anomaly flows}, \textit{Surveys in differential geometry 2018. {D}ifferential
              geometry, {C}alabi-{Y}au theory, and general relativity}, \textbf{23} (2020), pp. 89--103.

    \bibitem{FV2015} \textsc{\small A. Fino, L. Vezzoni}, \textup{\footnotesize Special Hermitian metrics on compact solvmanifolds}, \textit{J. Geom. Phys.}, \textbf{91} (2015), pp. 40--53.

    \bibitem{FLS2024} \textsc{\small E. Fusi, R. Lafuente, J. Stanfield}, \textup{\footnotesize The homogeneous generalized Ricci flow}, \textit{preprint arXiv:2404.15749}.
    
    \bibitem{LPV2020} \textsc{\small R. A. Lafuente, M. Pujia, L. Vezzoni}, \textup{\footnotesize Hermitian curvature flow on unimodular Lie groups and static invariant metrics}, \textit{Trans. Amer. Math. Soc.}, \textbf{373} (2020), pp. 3967--3993.

    \bibitem{L2012} \textsc{\small J. Lauret}, \textup{\footnotesize Convergence of homogeneous manifolds}, \textit{J. Lond. Math. Soc. (2)}, \textbf{86} (2012), pp. 701--727.

    \bibitem{L2015} \textsc{\small J. Lauret}, \textup{\footnotesize Curvature flows for almost-hermitian {L}ie groups}, \textit{Trans. Amer. Math. Soc.}, \textbf{367} (2015), pp. 7453--7480.

    \bibitem{L2013} \textsc{\small J. Lauret}, \textup{\footnotesize Ricci flow of homogeneous manifolds}, \textit{Math. Z.}, \textbf{274} (2013), pp. 373--403.

    \bibitem{L2011} \textsc{\small J. Lauret}, \textup{\footnotesize The {R}icci flow for simply connected nilmanifolds}, \textit{Comm. Anal. Geom.}, \textbf{19} (2011), pp. 831--854.

    \bibitem{L1993} \textsc{\small S. \L ojasiewicz}, \textup{\footnotesize Sur la g\'{e}om\'{e}trie semi- et sous-analytique}, \textit{Ann. Inst. Fourier (Grenoble)}, \textbf{43} (1993), pp. 1575--1595.
    
    \bibitem{PPZ2019} \textsc{\small D. H. Phong, S. Picard, X. Zhang}, \textup{\footnotesize A flow of conformally balanced metrics with {K}\"{a}hler fixed points}, \textit{Math. Ann.}, \textbf{374} (2019), pp. 2005--2040.

    \bibitem{PPZ2018} \textsc{\small D. H. Phong, S. Picard, X. Zhang}, \textup{\footnotesize Anomaly flows}, \textit{Comm. Anal. Geom.}, \textbf{26} (2018), pp. 955--1008.

    \bibitem{PPZ2018Strominger} \textsc{\small D. H. Phong, S. Picard, X. Zhang}, \textup{\footnotesize Geometric flows and {S}trominger systems}, \textit{Math. Z.}, \textbf{288} (2018), pp. 101--113.

    \bibitem{P2021} \textsc{\small M. Pujia}, \textup{\footnotesize Positive {H}ermitian curvature flow on complex 2-step nilpotent {L}ie groups}, \textit{Manuscripta Math.}, \textbf{166} (2021), pp. 237--249.

    \bibitem{R2009} \textsc{\small S. Rollenske}, \textup{\footnotesize Geometry of nilmanifolds with left-invariant complex structure and deformations in the large}, \textit{Proc. Lond. Math. Soc.}, \textbf{99} (2009), pp. 425–460.
    
    \bibitem{S2021} \textsc{\small J. Stanfield}, \textup{\footnotesize Positive {H}ermitian curvature flow on nilpotent and almost-abelian complex {L}ie groups}, \textit{Ann. Global Anal. Geom.}, \textbf{60} (2021), pp. 401--429.

    \bibitem{ST2011} \textsc{\small J. Streets, G. Tian}, \textup{\footnotesize Hermitian curvature flow}, \textit{J. Eur. Math. Soc. (JEMS)}, \textbf{13} (2011), pp. 601--634.

    \bibitem{ST2010} \textsc{\small J. Streets, G. Tian}, \textup{\footnotesize A parabolic flow of pluriclosed metrics}, \textit{Int. Math. Res. Not. IMRN}, \textbf{16} (2010), pp. 3101–-3133.

    \bibitem{U2019} \textsc{\small Y. Ustinovskiy}, \textup{\footnotesize The {H}ermitian curvature flow on manifolds with non-negative {G}riffiths curvature}, \textit{Amer. J. Math.}, \textbf{141} (2019), pp. 1751--1775.
    
\end{thebibliography}
\end{document}